\documentclass[11pt]{amsart}
\usepackage[utf8]{inputenc}
\usepackage{amsmath,amssymb,amsthm}
\usepackage{enumerate}
\usepackage[margin=2.7cm]{geometry}

\newtheorem{theorem}{Theorem}[section]
\numberwithin{equation}{section}
\newtheorem{lemma}[theorem]{Lemma}
\newtheorem{proposition}[theorem]{Proposition}
\newtheorem{corollary}[theorem]{Corollary}
\theoremstyle{definition}
\newtheorem{definition}[theorem]{Definition}
\newtheorem{example}[theorem]{Example}
\newtheorem*{definition*}{Definition}
\theoremstyle{remark}
\newtheorem{remark}[theorem]{Remark}

\newcommand{\Eh}{\hat{E}}
\newcommand{\Orth}{\operatorname{Orth}}

\begin{document}

\title{Pre-$d$-algebras, pre-$f$-algebras and their completions}

\author{Bahri Turan}
\address{Department of Mathematics, Faculty of Science, Gazi University, Ankara, T\"urkiye}
\email{bturan@gazi.edu.tr}
\urladdr{0000-0001-5256-1115}

\author{H\"uma G\"urk\"ok}
\address{H\"uma G\"urk\"ok}
\email{gurkokhuma@gmail.com}
\urladdr{0000-0002-3895-5115}

\subjclass[2020]{46A40, 06F25, 47B65}
\keywords{Pre-Riesz space, pre-$f$-algebra, pre-$d$-algebra, $f$-completion, Riesz completion, van Haandel extension, Riesz bimorphism, $d$-algebra, $f$-algebra, orthomorphism, disjointness}

\begin{abstract}
We define pre-$d$- and pre-$f$-algebras on Archimedean pre-Riesz spaces by conditions on the order and the multiplication only, and show that a pre-Riesz algebra is a pre-$f$-algebra if and only if its multiplication extends to an $f$-algebra multiplication on its Riesz completion. On order dense subspaces of $C(K)$, $K$ a compact Hausdorff space, in which disjoint positive elements separate the points of $K$, a pre-Riesz algebra is a pre-$f$-algebra if and only if multiplication by every positive element is band preserving.
\end{abstract}

\maketitle

\section{Introduction}

Pre-Riesz spaces are the partially ordered vector spaces that can be embedded as order dense subspaces of Riesz spaces; every directed Archimedean partially ordered vector space is of this kind. Van Haandel \cite{vH93} showed that a pre-Riesz space $E$ has a smallest such Riesz space, the Riesz completion $\hat E$, which is generated by $E$ as a Riesz space and is unique up to Riesz isomorphism. The Riesz completion makes it possible to transfer notions from Riesz spaces to pre-Riesz spaces and to compare them with their lattice counterparts. Many results on pre-Riesz spaces have been obtained in this way: a statement about $E$ is reformulated in $\hat E$, proved there by Riesz space methods, and then restricted to $E$ \cite{vGK06,vGK08,KvG19,KvGZ20,KK22,KSvG21}. In particular, disjointness, bands and ideals in pre-Riesz spaces were studied by Kalauch and van Gaans \cite{vGK06,vGK08,KvG19}, and disjointness preserving operators and their inverses in \cite{KvGZ20,KSvG21}. For operators, van Haandel's extension theorem states that a linear map between pre-Riesz spaces extends to a Riesz homomorphism between the Riesz completions if and only if it is a Riesz* homomorphism; Kalauch and Kusraeva extended this to Riesz $n$-morphisms \cite{KK22}.

In this paper we carry this approach over to algebraic structures. Let $A$ be an Archimedean pre-Riesz space with an associative bilinear multiplication under which products of positive elements are positive. We ask when this multiplication extends to the Riesz completion $\hat A$ in such a way that $\hat A$ becomes a $d$-algebra or an $f$-algebra, the classical lattice-ordered algebras studied by de Pagter, Huijsmans, Bernau and others \cite{BH90,HdP84,dP81}. Our aim is to answer this by conditions on $A$ alone, expressed through its order and its multiplication. This leads to the notions of pre-$d$-algebras and pre-$f$-algebras. The main tools are the extension theorem for Riesz $n$-morphisms of \cite{KK22} and the theory of $d$-algebras of Bernau and Huijsmans \cite{BH90}.

The paper is organised as follows. Section~\ref{sec:prelim} collects preliminaries on pre-Riesz spaces, Riesz* homomorphisms and lattice-ordered algebras. In Section~\ref{sec:pref} we introduce pre-$d$-algebras and prove that the multiplication of a pre-$d$-algebra $A$ extends uniquely to a $d$-algebra multiplication on $\hat A$, which gives the $d$-completion of $A$; we then introduce pre-$f$-algebras and show that this extension is an $f$-algebra multiplication if and only if $A$ is a pre-$f$-algebra. Finally we compare pre-$f$-algebras with local pre-Riesz algebras, which are defined by the condition in the definition of an $f$-algebra with lattice disjointness replaced by the disjointness of pre-Riesz spaces; locality is necessary but not sufficient for the pre-$f$ property. In Section~\ref{sec:f} we give further sufficient conditions for the pre-$f$ property, namely the existence of a unit or semiprimeness together with pervasiveness; we also show that a pre-$d$-algebra is commutative if and only if its squares are positive, and we characterise the $f$-algebra property of the $d$-completion by orthomorphisms. As an application of this characterisation we prove the main theorem: on order dense subspaces of $C(K)$, $K$ a compact Hausdorff space, in which disjoint positive elements separate the points of $K$, for instance on the spaces $C^k[0,1]$ $(k\in\mathbb N)$, locality and the pre-$f$ property coincide. Examples illustrating the results and the role of their hypotheses are given after the corresponding results.

\section{Preliminaries}\label{sec:prelim}

Let $E$ be a partially ordered vector space with positive cone $E_+$. For $M\subseteq E$ we denote by
\[
M^u=\{x\in E: x\ge m\text{ for all }m\in M\},\qquad M^l=\{x\in E: x\le m\text{ for all }m\in M\}
\]
the sets of upper and lower bounds of $M$, and we write $M^{ul}=(M^u)^l$; for $M=\{x,y\}$ we write $\{x,y\}^u$, $\{x,y\}^l$ and $\{x,y\}^{ul}$. The space $E$ is \emph{directed} if $\{x,y\}^u\neq\emptyset$ for all $x,y\in E$, and \emph{Archimedean} if $nx\le y$ for all $n\in\mathbb N$ implies $x\le0$ \cite[Definitions~1.1.6 and~1.1.9]{KvG19}. Let $E_1,E_2$ be partially ordered vector spaces. A linear map $T\colon E_1\to E_2$ is \emph{positive} if $Tx\ge0$ for all $x\in(E_1)_+$, and \emph{bipositive} if, for $x\in E_1$, $Tx\ge0$ holds if and only if $x\ge0$. A linear subspace $D$ of a partially ordered vector space $F$ is \emph{majorizing} if every element of $F$ has an upper bound in $D$, and \emph{order dense} if $y=\inf\{d\in D: d\ge y\}$ for every $y\in F$ \cite[Definition~1.6.1]{KvG19}; every order dense subspace is majorizing.

A partially ordered vector space $E$ is a \emph{pre-Riesz space} if for all $x,y,z\in E$ the inclusion $\{x+z,y+z\}^u\subseteq\{x,y\}^u$ implies $z\in E_+$. Every pre-Riesz space is directed, and every directed Archimedean partially ordered vector space is pre-Riesz \cite[Proposition~2.2.3]{KvG19}. A \emph{vector lattice cover} of $E$ is a pair $(F,j)$, where $F$ is a Riesz space and $j\colon E\to F$ is a bipositive linear map such that $j(E)$ is order dense in $F$; it is a \emph{Riesz completion} of $E$ if, in addition, no proper Riesz subspace of $F$ contains $j(E)$ \cite[Definition~2.4.6]{KvG19}. A Riesz completion $(\Eh,i)$ of $E$ exists exactly when $E$ is pre-Riesz, and it is unique up to Riesz isomorphism \cite[Theorem~2.4.5]{KvG19}. If $E$ is Archimedean, then so is every vector lattice cover of $E$ \cite[Proposition~2.4.12]{KvG19}.

For a compact Hausdorff space $K$, $C(K)$ denotes the Riesz space of real continuous functions on $K$ with the pointwise order, $\|x\|_\infty=\sup_{q\in K}|x(q)|$, and $\mathbf 1$ denotes the constant function with value $1$; for a compact interval $[s,t]$ we write $C[s,t]$. For $k\in\mathbb N$, $C^k[0,1]$ is the subspace of $C[0,1]$ consisting of the $k$ times continuously differentiable functions, with one-sided derivatives at $0$ and $1$. A function $x\in C[0,1]$ is \emph{piecewise $C^k$} if there is a partition $0=t_0<t_1<\dots<t_n=1$ such that the restriction of $x$ to each interval $[t_{j-1},t_j]$ is $k$ times continuously differentiable.

Let $E$ be a pre-Riesz space with Riesz completion $(\Eh,i)$. Since $\Eh$ is the Riesz subspace generated by $i(E)$, every element of $\Eh$ has the form
\[
\bigvee_{j=1}^{m} i(x_j)\;-\;\bigvee_{k=1}^{n} i(y_k),\qquad x_j,y_k\in E.
\]
The space $E$ is called \emph{pervasive} if for every $y\in\Eh_+$ with $y\neq 0$ there is $x\in E$ with $0<i(x)\le y$ \cite[Definition~2.8.1]{KvG19}. The spaces $C^k[0,1]$ with the pointwise order are pervasive, and their Riesz completions are the spaces of piecewise $C^k$ functions \cite[Example~2.8.2]{KvG19}.

Let $E_1,E_2$ be directed partially ordered vector spaces. A linear map $T\colon E_1\to E_2$ is a \emph{Riesz homomorphism} if $\bigl(T(\{x,y\}^u)\bigr)^l=\bigl(T(\{x,y\})\bigr)^{ul}$ for all $x,y\in E_1$ \cite[Definition~2.3.15]{KvG19}, and a \emph{Riesz* homomorphism} if $T(M^{ul})\subseteq\bigl(T(M)\bigr)^{ul}$ for every nonempty finite $M\subseteq E_1$ \cite[Definition~2.3.1]{KvG19}. Every Riesz homomorphism between pre-Riesz spaces is a Riesz* homomorphism \cite[Theorem~2.3.19]{KvG19}. If $E_1$ and $E_2$ are Riesz spaces, both notions coincide with that of a lattice homomorphism. A linear map $T\colon E_1\to E_2$ is a \emph{complete Riesz homomorphism} if $\inf M=0$ implies $\inf T(M)=0$ for every nonempty $M\subseteq E_1$; between Riesz spaces these are exactly the order continuous lattice homomorphisms, where a positive linear map $T\colon E_1\to E_2$ between Riesz spaces $E_1,E_2$ is \emph{order continuous} if $\inf T(D)=0$ for every downward directed set $D\subseteq E_1$ with $\inf D=0$.

\begin{theorem}[{\cite{vH93}, \cite[Theorems~2.3.12 and~2.4.11]{KvG19}}]\label{thm:vH}
Let $E_1,E_2$ be pre-Riesz spaces with Riesz completions $(\hat E_1,i_1)$, $(\hat E_2,i_2)$. A linear map $T\colon E_1\to E_2$ is a Riesz* homomorphism if and only if there is a Riesz homomorphism $\hat T\colon\hat E_1\to\hat E_2$ with $\hat T\circ i_1=i_2\circ T$; in this case $\hat T$ is unique.
\end{theorem}

The map $\hat T$ is called the \emph{van Haandel extension} of $T$. The uniqueness part rests on the following elementary observation, which we use repeatedly.

\begin{lemma}\label{lem:uniq}
Let $E$ be a pre-Riesz space with Riesz completion $(\Eh,i)$, let $F$ be a Riesz space, and let $S,T\colon\Eh\to F$ be lattice homomorphisms. Then $\{x\in\Eh: Sx=Tx\}$ is a Riesz subspace of $\Eh$. In particular, if $S\circ i=T\circ i$, then $S=T$.
\end{lemma}

\begin{proof}
Both maps are linear and preserve finite suprema, so the set where they agree is a Riesz subspace. If it contains $i(E)$, it is all of $\Eh$, since $\Eh$ is the Riesz subspace generated by $i(E)$.
\end{proof}

Let $E_1,\dots,E_n,F$ be partially ordered vector spaces. An $n$-linear map $T\colon E_1\times\dots\times E_n\to F$ is a \emph{Riesz $n$-morphism} if, for each $j\in\{1,\dots,n\}$ and each choice of positive entries $x_k\in(E_k)_+$ $(k\ne j)$, the map $E_j\to F$, $x\mapsto T(x_1,\dots,x_{j-1},x,x_{j+1},\dots,x_n)$, is a Riesz homomorphism \cite[Definition~3.1]{KK22}. A Riesz $2$-morphism is called a \emph{Riesz bimorphism}.

\begin{theorem}[{\cite[Theorem~3.7 and Corollary~3.8]{KK22}}]\label{thm:KK}
Let $E_1,\dots,E_n,F$ be pre-Riesz spaces with Riesz completions $(\hat E_1,i_1),\dots,(\hat E_n,i_n),(\hat F,i_F)$, and let $F$ be Archimedean. Every Riesz $n$-morphism $T\colon E_1\times\dots\times E_n\to F$ has a unique extension to a Riesz $n$-morphism $\hat T\colon\hat E_1\times\dots\times\hat E_n\to\hat F$, i.e.\ $\hat T(i_1(x_1),\dots,i_n(x_n))=i_F(T(x_1,\dots,x_n))$ for all $x_k\in E_k$. If $E_1=\dots=E_n$ and $T$ is symmetric, then $\hat T$ is symmetric.
\end{theorem}

We also use the following facts. Let $E$ be a pre-Riesz space with Riesz completion $(\Eh,i)$. If $F$ is a Riesz space, then a linear map $T\colon E\to F$ is a Riesz homomorphism if and only if
\[
T(x)\vee T(y)=\inf\{T(z): z\in E,\ z\ge x,y\}\qquad(x,y\in E)
\]
\cite[Proposition~2.3.26(ii)]{KvG19}, \cite[equation~(3)]{KK22}. If $F$ is a partially ordered vector space, $F_0$ is a directed linear subspace of $F$, and $S\colon E\to F$ is a Riesz homomorphism with $S(E)\subseteq F_0$, then the map $S_0\colon E\to F_0$, $S_0(x)=S(x)$, is a Riesz homomorphism \cite[Lemma~2.8]{KK22}. If $F$ is a partially ordered vector space and $j\colon E\to F$ is a bipositive linear map such that $j(E)$ is order dense in $F$, then $j$ is a complete Riesz homomorphism, hence a Riesz homomorphism \cite[Proposition~2.3.27, Theorem~2.3.19]{KvG19}; in particular $i\colon E\to\Eh$ is a Riesz homomorphism. The composition of two Riesz homomorphisms need not be a Riesz homomorphism \cite[Example~2.3.22]{KvG19}.

Let $E$ be a pre-Riesz space with Riesz completion $(\Eh,i)$. Elements $x,y\in E$ are \emph{disjoint}, $x\perp y$, if $\{x+y,-x-y\}^u=\{x-y,-x+y\}^u$; then $x\perp y$ in $E$ if and only if $i(x)\perp i(y)$ in $\Eh$ \cite[Proposition~4.1.4]{KvG19}. For $M\subseteq E$ let $M^d=\{x\in E: x\perp m\text{ for all }m\in M\}$; a set $B\subseteq E$ is a \emph{band} if $B=B^{dd}$ \cite[Definitions~4.1.1 and~4.2.1]{KvG19}. A linear operator $T\colon E\to E$ is \emph{disjointness preserving} if $x\perp y$ implies $Tx\perp Ty$ (\emph{on $E_+$} if this holds for $x,y\in E_+$), \emph{band preserving} if $T(B)\subseteq B$ for every band $B$ of $E$, and \emph{local} if $x\perp y$ implies $Tx\perp y$. Every local operator is disjointness preserving, and an operator is local if and only if it is band preserving \cite[Definition~5.1.1 and Proposition~5.1.2]{KvG19}. A set $M\subseteq E$ is \emph{order bounded} if it is contained in an order interval $[x,y]=\{z\in E: x\le z\le y\}$ with $x,y\in E$; a linear operator $T\colon E\to E$ is \emph{order bounded} if it maps order bounded sets to order bounded sets, and an \emph{orthomorphism} if it is order bounded and band preserving. We denote by $\Orth(E)$ the set of orthomorphisms of $E$ and by $I$ the identity operator on $E$. In Riesz spaces these notions are the classical ones \cite[Definition~1.4.13]{KvG19}, and for an Archimedean Riesz space $F$ every positive orthomorphism of $F$ is an order continuous lattice homomorphism \cite[Section~2]{HdP84}, \cite{dP81}, \cite[Chapter~20]{Za83}.

An \emph{$l$-algebra} is a Riesz space $A$ with an associative bilinear multiplication such that $A_+A_+\subseteq A_+$. Following \cite[Definition~1.1]{BH90}, $A$ is
\begin{enumerate}[(i)]
\item an \emph{almost $f$-algebra} if $a\wedge b=0$ implies $ab=0$ for all $a,b\in A$;
\item a \emph{$d$-algebra} if $c(a\vee b)=ca\vee cb$ and $(a\vee b)c=ac\vee bc$ for all $a,b\in A$ and $c\in A_+$, i.e.\ left and right multiplication by positive elements are lattice homomorphisms;
\item an \emph{$f$-algebra} if $a\wedge b=0$ implies $ca\wedge b=ac\wedge b=0$ for all $a,b\in A$ and $c\in A_+$, i.e.\ left and right multiplication by any element are orthomorphisms.
\end{enumerate}
For the theory of $f$-algebras and orthomorphisms we refer to de Pagter's thesis \cite{dP81}, to \cite[Section~2]{HdP84} and to \cite[Chapter~20]{Za83}; in particular, every Archimedean $f$-algebra is commutative. An element $a\in A$ is \emph{nilpotent} if $a^n=0$ for some $n\in\mathbb N$; $N(A)$ denotes the set of nilpotent elements of $A$, and $A$ is \emph{semiprime} if $N(A)=\{0\}$.

\begin{theorem}\label{thm:BH}
Let $A$ be an $l$-algebra.
\begin{enumerate}[(1)]
\item $A$ is a $d$-algebra if and only if $(ab)^+=a^+b^++a^-b^-$ and $(ab)^-=a^+b^-+a^-b^+$ for all $a,b\in A$ \cite[Proposition~1.2]{BH90}.
\item Every semiprime $d$-algebra is an $f$-algebra, and every $d$-algebra with a unit $e>0$ is an $f$-algebra \cite[Theorem~1.4]{BH90}.
\item If $A$ is an Archimedean $d$-algebra, then $A$ is commutative if and only if $a^2\ge 0$ for all $a\in A$, if and only if $A$ is an almost $f$-algebra \cite[Theorem~4.5]{BH90}.
\item If $A$ is an Archimedean $d$-algebra, then $N(A)$ is an $l$-ideal of $A$, i.e.\ an order ideal that is also a two-sided ring ideal \cite[Section~5]{BH90}.
\end{enumerate}
\end{theorem}

The following lemma is used to recognise multiplication operators on a pre-Riesz space as Riesz homomorphisms.

\begin{lemma}\label{lem:RH}
Let $E$ be a pre-Riesz space, let $(F,j)$ be a vector lattice cover of $E$, and let $S\colon F\to F$ be an order continuous lattice homomorphism with $S(j(E))\subseteq j(E)$. Then $T=j^{-1}\circ S\circ j\colon E\to E$ is a Riesz homomorphism.
\end{lemma}

\begin{proof}
Let $x,y\in E$. By bipositivity, $j(z)\ge j(x)\vee j(y)$ if and only if $z\ge x,y$, so order density gives $j(x)\vee j(y)=\inf\{j(z): z\in E,\ z\ge x,y\}$. An order continuous lattice homomorphism preserves every existing infimum: it preserves finite infima, and the finite infima of a set form a downward directed set with the same infimum. Hence
\[
Sj(x)\vee Sj(y)=S\bigl(j(x)\vee j(y)\bigr)=\inf\{Sj(z): z\in E,\ z\ge x,y\}.
\]
By \cite[equation~(3)]{KK22}, $S\circ j\colon E\to F$ is a Riesz homomorphism. Its image lies in the directed subspace $j(E)$ of $F$, so $S\circ j$ is a Riesz homomorphism $E\to j(E)$ by \cite[Lemma~2.8]{KK22}. Since $j$ is an order isomorphism of $E$ onto $j(E)$, $T=j^{-1}\circ S\circ j$ is a Riesz homomorphism $E\to E$.
\end{proof}

Let $K$ be a compact Hausdorff space. A subspace $E\subseteq C(K)$ \emph{strongly separates the points} of $K$ if for all $p\ne q$ in $K$ there is $x\in E$ with $x(q)=0$ and $x(p)=1$. If $E$ is order dense in $C(K)$, its Riesz completion is the Riesz subspace of $C(K)$ generated by $E$, with $i$ the inclusion map, and disjointness in $E$ is pointwise: for $x,y\in E$, $x\perp y$ if and only if $x(q)y(q)=0$ for all $q\in K$.

\begin{theorem}[{\cite[Theorem~5.1.14]{KvG19}}]\label{thm:vI}
Let $K$ be a compact Hausdorff space, let $E\subseteq C(K)$ be order dense and strongly separate the points of $K$, and let $T\colon E\to E$ be a Riesz* homomorphism. Then there are $w\in C(K)$, $w\ge0$, and $\alpha\colon K\to K$ such that $\alpha$ is continuous on the open set $\{w>0\}:=\{q\in K: w(q)>0\}$ and $(Tx)(q)=w(q)\,x(\alpha(q))$ for all $x\in E$ and $q\in K$.
\end{theorem}

In this situation $Sx=w\cdot(x\circ\alpha)$ $(x\in C(K))$ defines a lattice homomorphism $S\colon C(K)\to C(K)$; continuity of $Sx$ at points with $w(q)=0$ follows from $|Sx|\le w\|x\|_\infty$. Since $S$ maps $E$ into $E$, it maps $\Eh$ into $\Eh$, and by uniqueness its restriction to $\Eh$ is the van Haandel extension $\hat T$.

\section{The $d$-completion and the $f$-completion}\label{sec:pref}

In this section we introduce pre-$d$-algebras and pre-$f$-algebras. Our aim is to extend the multiplication of such an algebra to an $f$-algebra multiplication on the Riesz completion by means of conditions on the algebra itself. We first isolate a condition under which the multiplication extends to a $d$-algebra multiplication on the Riesz completion. In a $d$-algebra, multiplication by a positive element is a lattice homomorphism, $c(a\vee b)=ca\vee cb$; in a pre-Riesz space we replace $a\vee b$ by the set $\{a,b\}^u$ of upper bounds. All algebras considered are of the following kind.

\begin{definition}\label{def:alg}
An Archimedean pre-Riesz space $A$ is a \emph{pre-Riesz algebra} if it carries an associative bilinear multiplication with $A_+A_+\subseteq A_+$. For $c\in A$ we write $L_c,R_c\colon A\to A$ for the operators $L_c(a)=ca$ and $R_c(a)=ac$ $(a\in A)$, and for $M\subseteq A$ we write $cM=\{cb: b\in M\}$ and $Mc=\{bc: b\in M\}$. The Riesz completion of $A$ is denoted by $(\hat A,i)$.
\end{definition}

\begin{definition}\label{def:pred}
A pre-Riesz algebra $A$ is a \emph{pre-$d$-algebra} if for all $a,b\in A$ and $c\in A_+$,
\[
\bigl(c(\{a,b\}^u)\bigr)^l=\{ca,cb\}^{ul}\qquad\text{and}\qquad\bigl((\{a,b\}^u)c\bigr)^l=\{ac,bc\}^{ul}.
\]
\end{definition}

By the definition of a Riesz homomorphism recalled in Section~\ref{sec:prelim}, this says exactly that $L_c$ and $R_c$ are Riesz homomorphisms $A\to A$ for every $c\in A_+$; equivalently, the multiplication is a Riesz bimorphism $A\times A\to A$. If $A$ is a Riesz space, then $c(\{a,b\}^u)$ has the least element $c(a\vee b)$, so $\bigl(c(\{a,b\}^u)\bigr)^l=\{c(a\vee b)\}^l$ and $\{ca,cb\}^{ul}=\{ca\vee cb\}^l$; thus Definition~\ref{def:pred} becomes $c(a\vee b)=ca\vee cb$ and $(a\vee b)c=ac\vee bc$, the definition of a $d$-algebra in \cite{BH90}.

\begin{definition}\label{def:rhof}
Let $A$ be a pre-Riesz algebra. A \emph{$d$-completion} \textup{(}resp.\ an \emph{$f$-completion}\textup{)} of $A$ is a pair $(\hat A,i)$, where $\hat A$ is a $d$-algebra \textup{(}resp.\ an $f$-algebra\textup{)} and $i\colon A\to\hat A$ is a bipositive linear map such that
\begin{enumerate}[(i)]
\item $i(A)$ is order dense in $\hat A$;
\item $i(a)\cdot i(b)=i(ab)$ for all $a,b\in A$;
\item no proper Riesz subspace of $\hat A$ contains $i(A)$.
\end{enumerate}
\end{definition}

Conditions (i) and (iii) say that $(\hat A,i)$ is a Riesz completion of the pre-Riesz space $A$, and (ii) says that $i$ is an algebra homomorphism. Thus a $d$-completion \textup{(}$f$-completion\textup{)} is the Riesz completion together with a $d$-algebra \textup{(}$f$-algebra\textup{)} multiplication that extends the multiplication of $A$. Every $f$-algebra $B$ is a $d$-algebra (if $a\wedge b=0$ in $B$ and $c\in B_+$, then $ca\wedge b=0$, and applying the definition to the pair $b,ca$ gives $cb\wedge ca=0$; similarly on the right); hence every $f$-completion is a $d$-completion. Condition (iii) makes the completion unique.

\begin{proposition}\label{prop:funiq}
Let $A$ be a pre-Riesz algebra.
\begin{enumerate}[(a)]
\item If $(\hat A,i)$ is a $d$-completion of $A$, then the multiplication of $\hat A$ is the only $d$-algebra multiplication on the Riesz space $\hat A$ satisfying \textup{(ii)}.
\item If $(\hat A_1,i_1)$ and $(\hat A_2,i_2)$ are $d$-completions of $A$, then there is a unique Riesz isomorphism $\varphi\colon\hat A_1\to\hat A_2$ with $\varphi\circ i_1=i_2$, and $\varphi$ is an algebra isomorphism.
\end{enumerate}
In particular, both statements hold for $f$-completions.
\end{proposition}

\begin{proof}
(a) Let $\cdot$ and $\star$ be $d$-algebra multiplications on $\hat A$ satisfying (ii). For $c\in A_+$ the maps $\hat A\to\hat A$ given by $b\mapsto b\cdot i(c)$ and $b\mapsto b\star i(c)$ are lattice homomorphisms that agree on $i(A)$, both sending $i(a)$ to $i(ac)$ $(a\in A)$. By (iii) and Lemma~\ref{lem:uniq} they agree on $\hat A$. As $A$ is directed, $b_1\cdot b_2=b_1\star b_2$ for all $b_1\in\hat A$ and $b_2\in i(A)$. Now fix $b_1\in\hat A_+$. The maps $\hat A\to\hat A$ given by $b\mapsto b_1\cdot b$ and $b\mapsto b_1\star b$ are lattice homomorphisms that agree on $i(A)$, hence on $\hat A$. Writing $b_1=b_1^+-b_1^-$ gives $b_1\cdot b_2=b_1\star b_2$ for all $b_1,b_2\in\hat A$.

(b) By (i) and (iii), $(\hat A_1,i_1)$ and $(\hat A_2,i_2)$ are Riesz completions of $A$, so there is a unique Riesz isomorphism $\varphi\colon\hat A_1\to\hat A_2$ with $\varphi\circ i_1=i_2$ \cite[Section~2.4]{KvG19}. For $b_1,b_2\in\hat A_2$ define $b_1\star b_2=\varphi\bigl(\varphi^{-1}(b_1)\cdot\varphi^{-1}(b_2)\bigr)$, where $\cdot$ is the multiplication of $\hat A_1$. Since $\varphi$ is a Riesz isomorphism, $\star$ is a $d$-algebra multiplication on $\hat A_2$, and $i_2(a)\star i_2(b)=\varphi(i_1(a)\cdot i_1(b))=\varphi(i_1(ab))=i_2(ab)$. By (a), $\star$ is the multiplication of $\hat A_2$, i.e.\ $\varphi$ is multiplicative.
\end{proof}

By Proposition~\ref{prop:funiq} we may speak of \emph{the} $d$-completion and \emph{the} $f$-completion of $A$, and keep the notation $\hat A$ of the Riesz completion. We now show that every pre-$d$-algebra has a $d$-completion.

\begin{lemma}\label{lem:tri}
Let $A$ be a pre-Riesz algebra with Riesz completion $(\hat A,i)$, let $F$ be a Riesz space, and let $S,T\colon\hat A\times\hat A\times\hat A\to F$ be Riesz $3$-morphisms. If $S$ and $T$ agree on $i(A)\times i(A)\times i(A)$, then $S=T$.
\end{lemma}

\begin{proof}
Fix $b_2,b_3\in i(A_+)$. The maps $\hat A\to F$ given by $b\mapsto S(b,b_2,b_3)$ and $b\mapsto T(b,b_2,b_3)$ are lattice homomorphisms that agree on $i(A)$, so they agree on $\hat A$ by Lemma~\ref{lem:uniq}. Since $A$ is directed, $S$ and $T$ agree on $\hat A\times i(A)\times i(A)$. Next fix $b_1\in\hat A_+$ and $b_3\in i(A_+)$; the maps $\hat A\to F$ given by $b\mapsto S(b_1,b,b_3)$ and $b\mapsto T(b_1,b,b_3)$ are lattice homomorphisms agreeing on $i(A)$, hence on $\hat A$. Writing $b_1=b_1^+-b_1^-$ and $b_3$ as a difference of elements of $i(A_+)$, $S$ and $T$ agree on $\hat A\times\hat A\times i(A)$. The same argument, applied for $b_1,b_2\in\hat A_+$ to the maps $\hat A\to F$ given by $b\mapsto S(b_1,b_2,b)$ and $b\mapsto T(b_1,b_2,b)$, completes the proof.
\end{proof}

\begin{theorem}\label{thm:ext}
Every pre-$d$-algebra $A$ has a $d$-completion. More precisely, there is a unique bilinear multiplication on $\hat A$ that is a Riesz bimorphism and satisfies $i(a)\cdot i(b)=i(ab)$ for $a,b\in A$; with this multiplication $\hat A$ is an Archimedean $d$-algebra, and $(\hat A,i)$ is the $d$-completion of $A$. If $A$ is commutative, so is $\hat A$, and then $\hat A$ is an almost $f$-algebra.
\end{theorem}

\begin{proof}
The multiplication $m\colon A\times A\to A$ is a Riesz bimorphism and $A$ is Archimedean by Definition~\ref{def:alg}, so Theorem~\ref{thm:KK} gives a unique Riesz bimorphism $\hat m\colon\hat A\times\hat A\to\hat A$ with $\hat m(i(a),i(b))=i(ab)$. Write $b_1\cdot b_2=\hat m(b_1,b_2)$ for $b_1,b_2\in\hat A$. For $b_2\in\hat A_+$ the maps $\hat A\to\hat A$ given by $b\mapsto b\cdot b_2$ and $b\mapsto b_2\cdot b$ are lattice homomorphisms, and products of positive elements are positive.

It remains to prove associativity. For $b_1,b_2,b_3\in\hat A$ define $S(b_1,b_2,b_3)=(b_1\cdot b_2)\cdot b_3$ and $T(b_1,b_2,b_3)=b_1\cdot(b_2\cdot b_3)$. Then $S,T\colon\hat A\times\hat A\times\hat A\to\hat A$ are trilinear. For $b_2,b_3\in\hat A_+$ the map $\hat A\to\hat A$, $b\mapsto(b\cdot b_2)\cdot b_3$, is the composition of the lattice homomorphisms $b\mapsto b\cdot b_2$ and $b\mapsto b\cdot b_3$ of $\hat A$, hence a lattice homomorphism; the same holds in the other variables and for $T$. So $S$ and $T$ are Riesz $3$-morphisms. On $i(A)^3$ they agree by associativity in $A$, since $i(a)\cdot i(b)=i(ab)$. Lemma~\ref{lem:tri} gives $S=T$. Hence $\hat A$ is an $l$-algebra in which multiplication by positive elements is a lattice homomorphism, i.e.\ a $d$-algebra; it is Archimedean, since $A$ is \cite[Proposition~2.4.12]{KvG19}.

If $A$ is commutative, $m$ is symmetric and $\hat m$ is symmetric by Theorem~\ref{thm:KK}. Then $\hat A$ is a commutative Archimedean $d$-algebra, hence an almost $f$-algebra by Theorem~\ref{thm:BH}(3). Finally, $(\hat A,i)$ satisfies (i)--(iii) of Definition~\ref{def:rhof}, so it is a $d$-completion of $A$, and by Proposition~\ref{prop:funiq} it is the only one.
\end{proof}

\begin{remark}\label{rem:vH}
For $c\in A_+$, $L_c$ and $R_c$ are Riesz homomorphisms, hence Riesz* homomorphisms, and by uniqueness in Theorem~\ref{thm:vH} the maps $\hat A\to\hat A$ given by $b\mapsto i(c)\cdot b$ and $b\mapsto b\cdot i(c)$ are their van Haandel extensions:
\[
b\cdot i(c)=\hat R_c(b),\qquad i(c)\cdot b=\hat L_c(b)\qquad(b\in\hat A,\ c\in A_+).
\]
For $b=\bigvee_j i(a_j)-\bigvee_k i(a_k')$ with $a_j,a_k'\in A$ one has $b\cdot i(c)=\bigvee_j i(a_jc)-\bigvee_k i(a_k'c)$.
\end{remark}

Conversely, if $\hat A$ carries a $d$-algebra multiplication extending that of $A$, then $L_c$ and $R_c$ $(c\in A_+)$ are Riesz* homomorphisms by Theorem~\ref{thm:vH}. The following example shows that they need not be Riesz homomorphisms.

\begin{example}\label{ex:star}
Let $A=\{a\in C[-1,1]: a(-1)+a(1)=2a(0)\}$, which is order dense in $C[-1,1]$ \cite[Example~2.3.22]{KvG19}, and define $a*b=b(0)\,a$ for $a,b\in A$.
\begin{enumerate}[(1)]
\item The multiplication is associative and positive, and $R_b=b(0)I$ is a Riesz homomorphism for $b\in A_+$.
\item For $a\in A_+$, $a\ne0$, the operator $L_a(b)=b(0)a$ is a Riesz* homomorphism, being the restriction of the lattice homomorphism $C[-1,1]\to C[-1,1]$, $x\mapsto x(0)a$, \cite[Corollary~2.3.7]{KvG19}. It is not a Riesz homomorphism. Let $b_1,b_2\in A$ be given by $b_1(s)=s$ and $b_2(s)=-s$ $(s\in[-1,1])$. Every $b\in\{b_1,b_2\}^u$ satisfies $b(\pm1)\ge1$ and hence $b(0)\ge1$, so $a\in\bigl(L_a(\{b_1,b_2\}^u)\bigr)^l$, whereas $\{L_ab_1,L_ab_2\}^{ul}=\{0\}^{ul}=-A_+$ does not contain $a$. So $A$ is not a pre-$d$-algebra in the sense of Definition~\ref{def:pred}.
\item Nevertheless, $b\cdot b'=b'(0)\,b$ $(b,b'\in\hat A)$ is a $d$-algebra multiplication on $\hat A$, the Riesz subspace of $C[-1,1]$ generated by $A$, and it extends the multiplication of $A$: point evaluation at $0$ is a lattice homomorphism of $\hat A$.
\end{enumerate}
This space is the Namioka space of \cite[Example~42]{KSvG21}.
\end{example}

By Theorem~\ref{thm:ext}, every pre-$d$-algebra has a $d$-completion. We now introduce the condition under which it is an $f$-algebra. In an Archimedean $f$-algebra $A$ one has $a^+\perp(ac)^-$ and $a^+\perp(ca)^-$ for all $a\in A$ and $c\in A_+$, since $(ac)^-=a^-c$ and $(ca)^-=ca^-$ belong to the band generated by $a^-$. In a pre-Riesz algebra $A$ the elements $i(a)^+$ and $i(ac)^-$ of $\hat A$ $(a\in A,\ c\in A_+)$ are represented by the sets $\{0,a\}^u$ and $\{0,-ac\}^u$ of upper bounds, and we formulate the condition by means of these sets.

\begin{definition}\label{def:pref}
A pre-$d$-algebra $A$ is a \emph{pre-$f$-algebra} if for all $a\in A$ and $c\in A_+$,
\[
\bigl(\{0,a\}^u\cup\{0,-ac\}^u\bigr)^l\subseteq -A_+\qquad\text{and}\qquad \bigl(\{0,a\}^u\cup\{0,-ca\}^u\bigr)^l\subseteq -A_+.
\]
\end{definition}

We call the condition in Definition~\ref{def:pref} the \emph{sign condition}. Like the pre-$d$ condition, it refers only to the order and the multiplication of $A$. If $A$ is a Riesz space, the sign condition says $a^+\wedge(ac)^-=a^+\wedge(ca)^-=0$; in a $d$-algebra $(ac)^-=a^-c$ and $(ca)^-=ca^-$ for $c\in A_+$ by Theorem~\ref{thm:BH}(1), so the condition becomes $a^+\wedge a^-c=a^+\wedge ca^-=0$, which is the definition of an $f$-algebra applied to the disjoint pair $(a^+,a^-)$. Hence an Archimedean $l$-algebra is a pre-$f$-algebra if and only if it is an $f$-algebra.

The following example shows that a pre-$d$-algebra need not be a pre-$f$-algebra and need not have an $f$-completion.

\begin{example}\label{ex:nofcomp}
Let $A$ be the space of real trigonometric polynomials on the circle $S^1=\mathbb R/2\pi\mathbb Z$, i.e.\ the linear combinations of $\mathbf 1$, $\cos k\theta$ and $\sin k\theta$ $(k\in\mathbb N)$, with the pointwise order. Fix $\theta_0\in S^1$, and define $c*a=c(\theta_0)\,a(\theta_0)\mathbf 1$ for $a,c\in A$. As a subspace of $C(S^1)$, $A$ is Archimedean, and it is directed, since $\mathbf 1\in A$ and $a\le\|a\|_\infty\mathbf 1$ for $a\in A$; hence $A$ is a pre-Riesz space \cite[Proposition~2.2.3]{KvG19}. The multiplication $*$ is associative, commutative and positive, so $A$ is a pre-Riesz algebra. By the Weierstrass theorem for trigonometric polynomials, for every $x\in C(S^1)$ and $\varepsilon>0$ there is $a\in A$ with $x\le a\le x+\varepsilon\mathbf 1$; hence $A$ is order dense in $C(S^1)$, and its Riesz completion is the Riesz subspace $\hat A$ of $C(S^1)$ generated by $A$, with $i$ the inclusion.

We first show that $A$ is a pre-$d$-algebra. For $c\in A_+$ and $a\in A$ we have $L_c(a)=R_c(a)=c(\theta_0)\,a(\theta_0)\mathbf 1$. Let $a_1,a_2\in A$ and $M=\max(a_1(\theta_0),a_2(\theta_0))$. For $\varepsilon>0$, the Weierstrass theorem for trigonometric polynomials gives $b\in A$ with $\max(a_1,a_2)\le b\le\max(a_1,a_2)+\varepsilon\mathbf 1$, so $\inf\{b(\theta_0): b\in A,\ b\ge a_1,a_2\}=M$, and since $c(\theta_0)\ge0$,
\[
L_c(a_1)\vee L_c(a_2)=c(\theta_0)M\,\mathbf 1=\inf\{L_c(b): b\in A,\ b\ge a_1,a_2\}
\]
in $\hat A$. By the characterisation of Riesz homomorphisms recalled in Section~\ref{sec:prelim}, $L_c=R_c$ is a Riesz homomorphism $A\to\hat A$ with image in the directed subspace $A$, hence a Riesz homomorphism $A\to A$ \cite[Lemma~2.8]{KK22}.

Next we show that $A$ is not a pre-$f$-algebra. Let $a(\theta)=\tfrac12-\cos(\theta-\theta_0)$ and $c=\mathbf 1$. Then $a*c=a(\theta_0)\mathbf 1=-\tfrac12\mathbf 1$, so $\{0,-a*c\}^u=\{a'\in A: a'\ge\tfrac12\mathbf 1\}$. Put $b(\theta)=-\tfrac12\cos(\theta-\theta_0)$. Then $b\le\tfrac12\mathbf 1$ and $a-b=\tfrac12\bigl(1-\cos(\theta-\theta_0)\bigr)\ge0$. Hence $b$ is a lower bound of $\{0,a\}^u\cup\{0,-a*c\}^u$, but $b(\theta_0+\pi)=\tfrac12>0$, so $b\notin-A_+$.

Finally, we show that $A$ has no $f$-completion. Point evaluation at $\theta_0$ is a lattice homomorphism of $\hat A\subseteq C(S^1)$, so $b\cdot b'=b(\theta_0)b'(\theta_0)\mathbf 1$ $(b,b'\in\hat A)$ defines a Riesz bimorphism $\hat A\times\hat A\to\hat A$ extending $*$; by Theorem~\ref{thm:ext} it is the multiplication of the $d$-completion. It is not an $f$-algebra multiplication: for $b_1=(\cos(\theta-\theta_0))^+$ and $b_2=(\cos(\theta-\theta_0))^-$ in $\hat A$ we have $b_1\wedge b_2=0$, but $(b_1\cdot\mathbf 1)\wedge b_2=\mathbf 1\wedge b_2=b_2\ne0$. Since every $f$-completion is a $d$-completion and the $d$-completion is unique (Proposition~\ref{prop:funiq}), $A$ has no $f$-completion.
\end{example}

We now show that the pre-$f$-algebras are exactly the pre-Riesz algebras with an $f$-completion. First, an $f$-completion can exist only for pre-$d$-algebras.

\begin{proposition}\label{prop:rhof}
Let $A$ be a pre-Riesz algebra with an $f$-completion $(\hat A,i)$. Then $A$ is a commutative pre-$d$-algebra, and $(\hat A,i)$ is the $d$-completion of $A$ given by Theorem~\ref{thm:ext}.
\end{proposition}

\begin{proof}
Let $c\in A_+$. Since $A$ is Archimedean, so is $\hat A$ \cite[Proposition~2.4.12]{KvG19}. Right multiplication by $i(c)$ is a positive orthomorphism of the Archimedean Riesz space $\hat A$, hence an order continuous lattice homomorphism, and it maps $i(A)$ into $i(A)$. By Lemma~\ref{lem:RH}, applied to the cover $(\hat A,i)$, $R_c$ is a Riesz homomorphism $A\to A$; the same holds for $L_c$, so $A$ is a pre-$d$-algebra. Since $\hat A$ is an $f$-algebra, it is a $d$-algebra, so $(\hat A,i)$ is a $d$-completion of $A$, and by Proposition~\ref{prop:funiq}(b) it coincides with the one of Theorem~\ref{thm:ext}. Finally, $\hat A$ is an Archimedean $f$-algebra, hence commutative, and $i$ is injective and multiplicative, so $A$ is commutative.
\end{proof}

\begin{corollary}\label{cor:fiffd}
A pre-Riesz algebra $A$ has an $f$-completion if and only if $A$ is a pre-$d$-algebra and its $d$-completion is an $f$-algebra.
\end{corollary}

\begin{proof}
If $A$ has an $f$-completion, then $A$ is a pre-$d$-algebra and its $d$-completion is an $f$-algebra by Proposition~\ref{prop:rhof}. Conversely, if $A$ is a pre-$d$-algebra whose $d$-completion is an $f$-algebra, then this $d$-completion is an $f$-completion of $A$ by Definition~\ref{def:rhof}.
\end{proof}

By Corollary~\ref{cor:fiffd}, it remains to show that for a pre-$d$-algebra the condition of Definition~\ref{def:pref} is equivalent to the $f$-algebra property of the $d$-completion. The proof of the following theorem is based on a description of the positive and negative parts of elements of $\hat A$, which together with the Riesz bimorphism property of the multiplication reduces disjointness conditions on $\hat A$ to elements of $A$. In Theorem~\ref{thm:sign} and Corollary~\ref{cor:intr}, $A$ is a pre-$d$-algebra and $\hat A$ carries the multiplication of Theorem~\ref{thm:ext}.

\begin{theorem}\label{thm:sign}
Let $A$ be a pre-$d$-algebra. Then $\hat A$ is an $f$-algebra if and only if, for all $a\in A$ and $c\in A_+$,
\[
i(a)^+\perp i(ac)^-\qquad\text{and}\qquad i(a)^+\perp i(ca)^-.
\]
\end{theorem}

\begin{proof}
Suppose first that $\hat A$ is an $f$-algebra. Then $i(ac)^-=\hat R_c(i(a)^-)=i(a)^-\cdot i(c)$ lies in the band generated by $i(a)^-$, which is disjoint from $i(a)^+$; similarly for $ca$.

Conversely, suppose that $i(a)^+\perp i(ac)^-$ and $i(a)^+\perp i(ca)^-$ for all $a\in A$ and $c\in A_+$. Note that $\hat A$ is an $f$-algebra if and only if $b^+\perp b^-\cdot b'$ and $b^+\perp b'\cdot b^-$ for all $b\in\hat A$ and $b'\in\hat A_+$, because every pair $(b_1,b_2)$ of disjoint elements of $\hat A_+$ equals $(b^-,b^+)$ for $b=b_2-b_1$. Fix $b\in\hat A$. Since the map $\hat A\to\hat A$, $b'\mapsto b^-\cdot b'$, is a lattice homomorphism and $\{b^+\}^d$ is a band, the set $\{b'\in\hat A: b^-\cdot b'\perp b^+\}$ is a Riesz subspace of $\hat A$. It contains $i(A_+)$ if it contains $i(c)$ for every $c\in A_+$, and then it contains $i(A)$, since $A$ is directed, and hence equals $\hat A$, the Riesz subspace generated by $i(A)$. So it suffices to take $b'=i(c)$ with $c\in A_+$.

Now write $b=\bigvee_{j=1}^{m}i(a_j)-\bigvee_{k=1}^{n}i(a_k')$, where $a_j,a_k'\in A$, and put $a_{jk}=a_j-a_k'\in A$. Then $b^+=\bigvee_j\bigl(i(a_j)-\bigvee_k i(a_k')\bigr)^+$ and $\bigl(i(a_j)-\bigvee_k i(a_k')\bigr)^+=\bigwedge_k i(a_{jk})^+$; in the same way $b^-=\bigvee_k\bigwedge_j\bigl(i(a_k')-i(a_j)\bigr)^+$, and $\bigl(i(a_k')-i(a_j)\bigr)^+=i(a_{jk})^-$. Hence
\begin{equation}\label{eq:red}
b^+=\bigvee_{j=1}^{m}\bigwedge_{k=1}^{n} i(a_{jk})^+,\qquad b^-=\bigvee_{k=1}^{n}\bigwedge_{j=1}^{m} i(a_{jk})^-.
\end{equation}
Right multiplication by $i(c)$ is a lattice homomorphism, and finite suprema and infima distribute, so
\[
b^+\wedge b^-\,i(c)\ \le\ \bigvee_{j,k}\bigl(i(a_{jk})^+\wedge i(a_{jk})^-\,i(c)\bigr).
\]
For $a\in A$ we have $i(a)^-\cdot i(c)=\hat R_c(i(a)^-)=(\hat R_c\,i(a))^-=i(ac)^-$, so each term vanishes by hypothesis, applied to $a=a_{jk}$. Left multiplication is treated in the same way, using $\hat L_c$ and $ca$.
\end{proof}

The condition of Theorem~\ref{thm:sign} can be expressed in $A$ alone; it is the condition of Definition~\ref{def:pref}.

\begin{corollary}\label{cor:intr}
Let $A$ be a pre-$d$-algebra. Then $\hat A$ is an $f$-algebra if and only if, for all $a\in A$ and $c\in A_+$,
\[
\bigl(\{0,a\}^u\cup\{0,-ac\}^u\bigr)^l\subseteq -A_+\qquad\text{and}\qquad \bigl(\{0,a\}^u\cup\{0,-ca\}^u\bigr)^l\subseteq -A_+.
\]
\end{corollary}

\begin{proof}
By order density, $i(a)^+=\inf\{i(b): b\in\{0,a\}^u\}$ and $i(ac)^-=\inf\{i(b): b\in\{0,-ac\}^u\}$ in $\hat A$. Hence $b\in A$ is a lower bound of $\{0,a\}^u\cup\{0,-ac\}^u$ if and only if $i(b)\le b_0:=i(a)^+\wedge i(ac)^-$. If $b_0=0$, every such $b$ satisfies $b\le0$. Conversely, applying order density to $-b_0$ gives $b_0=\sup\{i(b): b\in A,\ i(b)\le b_0\}$. If $b_0>0$, not all these $b$ can lie in $-A_+$, since otherwise the supremum would be $\le0$; so there is a lower bound outside $-A_+$. The second condition is treated in the same way, and Theorem~\ref{thm:sign} completes the proof.
\end{proof}

Since $(\{0,a\}^u\cup\{0,-ac\}^u)^l=\{0,a\}^{ul}\cap\{0,-ac\}^{ul}$, the condition of Corollary~\ref{cor:intr} says that $\{0,a\}^{ul}\cap\{0,-ac\}^{ul}\subseteq-A_+$, where $\{0,a\}^{ul}$ and $\{0,-ac\}^{ul}$ are the sets of lower bounds that represent $i(a)^+$ and $i(ac)^-$ in van Haandel's construction of $\hat A$; this is the description of disjointness of positive elements of $\hat A$ in \cite[Section~6]{vGK06}.

\begin{theorem}\label{thm:fcomp}
For a pre-Riesz algebra $A$ the following are equivalent.
\begin{enumerate}[(i)]
\item $A$ is a pre-$f$-algebra.
\item $A$ has an $f$-completion.
\end{enumerate}
In this case the $f$-completion is unique up to an isomorphism of $f$-algebras, it coincides with the $d$-completion of Theorem~\ref{thm:ext}, and $A$ is commutative.
\end{theorem}

\begin{proof}
(i)$\Rightarrow$(ii). By Theorem~\ref{thm:ext}, $A$ has a $d$-completion $(\hat A,i)$. By Corollary~\ref{cor:intr} and the sign condition of Definition~\ref{def:pref}, $\hat A$ is an $f$-algebra, so $(\hat A,i)$ is an $f$-completion.

(ii)$\Rightarrow$(i). By Proposition~\ref{prop:rhof}, $A$ is a pre-$d$-algebra and its $f$-completion is its $d$-completion. As this is an $f$-algebra, Corollary~\ref{cor:intr} gives the sign condition.

The remaining statements follow from Propositions~\ref{prop:funiq} and~\ref{prop:rhof}.
\end{proof}

If $A$ is an Archimedean $f$-algebra, then $(A,i)$, with $i$ the identity map of $A$, is its $f$-completion, and the zero multiplication on any Archimedean pre-Riesz space is a pre-$f$-algebra structure. Further sufficient conditions are given in Section~\ref{sec:f}.

The following example gives a pre-Riesz algebra that is not a Riesz space and has an $f$-completion.

\begin{example}\label{ex:fcomp}
Let $A=P[0,1]$ be the space of real polynomials on $[0,1]$ with the pointwise order and the pointwise product. Then $A$ is order dense in $C[0,1]$ by the Weierstrass theorem, and its Riesz completion is the Riesz subspace $PP[0,1]$ of $C[0,1]$ consisting of the continuous piecewise polynomial functions, with $i$ the inclusion map \cite[Example~3.5]{vGK06}. A product of two continuous piecewise polynomial functions is again one, so $PP[0,1]$ is both a Riesz subspace and a subalgebra of the Archimedean $f$-algebra $C[0,1]$. Lattice operations and multiplication are inherited from $C[0,1]$, so $PP[0,1]$ is an $f$-algebra. In the terms of Theorem~\ref{thm:fcomp}: $A$ is a pre-$d$-algebra by Lemma~\ref{lem:RH}, since for $c\in A_+$ multiplication by $c$ is an order continuous lattice homomorphism $C[0,1]\to C[0,1]$ that maps $A$ into $A$, and it is a pre-$f$-algebra, since for $a\in A$ and $c\in A_+$ the functions $a$ and $ac$ have the same sign at every $t\in[0,1]$ with $c(t)>0$, so that $i(a)^+\perp i(ac)^-$; its $d$-completion is $PP[0,1]$ with the pointwise product, which is an $f$-algebra. Hence $(PP[0,1],i)$ is the $f$-completion of $A$, although $A$ is not a Riesz space.
\end{example}

The definition of an $f$-algebra uses the lattice structure only through disjointness, and pre-Riesz spaces carry the intrinsic disjointness of \cite{vGK06}. We compare pre-$f$-algebras with the algebras defined by the disjointness condition in the definition of an $f$-algebra.

\begin{definition}\label{def:local}
A pre-Riesz algebra $A$ is \emph{local} if for all $a,b\in A$ and $c\in A_+$,
\[
a\perp b\ \Longrightarrow\ ca\perp b\ \text{ and }\ ac\perp b .
\]
\end{definition}

\begin{lemma}\label{lem:orth}
A pre-Riesz algebra $A$ is local if and only if $L_c,R_c\in\Orth(A)$ for every $c\in A_+$.
\end{lemma}

\begin{proof}
For $c\in A_+$ the operators $L_c,R_c$ are positive, hence order bounded. The condition of Definition~\ref{def:local} says that $L_c$ and $R_c$ are local, and an operator on a pre-Riesz space is local if and only if it is band preserving \cite[Proposition~5.1.2]{KvG19}.
\end{proof}

If $A$ is a Riesz space, Definition~\ref{def:local} is the definition of an $f$-algebra \cite[Definition~1.1]{BH90}, \cite[Section~2]{HdP84}, since the disjointness of \cite{vGK06} coincides with the classical one. Lemma~\ref{lem:orth} is the analogue for pre-Riesz algebras of the description of $f$-algebras by orthomorphisms. Since local operators are disjointness preserving, in a local pre-Riesz algebra the operators $L_c,R_c$ $(c\in A_+)$ are positive and disjointness preserving.

\begin{proposition}\label{prop:pflocal}
Every pre-$f$-algebra is a local pre-Riesz algebra.
\end{proposition}

\begin{proof}
Let $A$ be a pre-$f$-algebra with $f$-completion $(\hat A,i)$ (Theorem~\ref{thm:fcomp}), and let $a,b\in A$ with $a\perp b$ and $c\in A_+$. Then $i(a)\perp i(b)$ \cite[Proposition~4.1.4]{KvG19}, so $i(ca)=i(c)\cdot i(a)\perp i(b)$ since $\hat A$ is an $f$-algebra, and hence $ca\perp b$; similarly $ac\perp b$.
\end{proof}

The converse fails. The algebra of Example~\ref{ex:nofcomp} is local, because it has no nonzero disjoint elements (disjointness in $A$ is pointwise, and a nonzero trigonometric polynomial has only finitely many zeros), but it is not a pre-$f$-algebra. In a Riesz space, positive disjointness preserving operators are lattice homomorphisms, so every $f$-algebra is a $d$-algebra. The following example shows that a local pre-Riesz algebra need not be a pre-$d$-algebra.

\begin{example}\label{ex:int}
Let $A=P[0,1]$ be the space of real polynomials on $[0,1]$ with the pointwise order, and define $a*b=\bigl(\int_0^1a\bigr)\bigl(\int_0^1b\bigr)\mathbf 1$ for $a,b\in A$. The multiplication is associative, commutative and positive. Since $A$ has no nonzero disjoint elements, $A$ is a local pre-Riesz algebra. For $c\in A_+$ with $\int_0^1c=1$ we have $L_c(a)=\bigl(\int_0^1a\bigr)\mathbf 1$. This is not a Riesz* homomorphism: $A$ strongly separates the points of $[0,1]$, so by Theorem~\ref{thm:vI} there would be $w\in C[0,1]$, $w\ge0$, and $\alpha\colon[0,1]\to[0,1]$ with $\int_0^1a=w(q)\,a(\alpha(q))$ for all $a\in A$ and $q\in[0,1]$; the choices $a=\mathbf 1$, $a(t)=t$ and $a(t)=t^2$ give $w=1$, $\alpha(q)=\tfrac12$ and $\tfrac13=\tfrac14$, a contradiction (compare \cite[Example~5.1.15]{KvG19}). Hence $A$ is a local pre-Riesz algebra but not a pre-$d$-algebra, and in particular not a pre-$f$-algebra.
\end{example}

Conversely, a pre-$d$-algebra need not be local, even if it is a Riesz space.

\begin{example}\label{ex:R2}
Let $A=\mathbb R^2$ with the coordinatewise order and $a\cdot b=(a_1b_1,a_1b_1)$ for $a=(a_1,a_2)$ and $b=(b_1,b_2)$ in $A$. This is a commutative Archimedean $d$-algebra that is not an $f$-algebra: for the unit vectors $e_1=(1,0)$ and $e_2=(0,1)$ we have $e_1\wedge e_2=0$ but $(e_1\cdot e_1)\wedge e_2=e_2$. It has no unit and $e_2$ is nilpotent. Here $\hat A=A$, so $A$ is a pre-$d$-algebra but not a local pre-Riesz algebra; compare \cite[Example~3.3]{BH90}.
\end{example}

The algebra of Example~\ref{ex:star}, which has a $d$-completion, is not local either: if $b\in A$ satisfies $b(0)\ne0$ and vanishes outside $[-1,-\tfrac9{10}]\cup[-\tfrac1{10},\tfrac1{10}]\cup[\tfrac9{10},1]$, and $a\in A\setminus\{0\}$ vanishes outside $(\tfrac15,\tfrac45)$, then $a\perp b$ but $\mathbf 1*b=b(0)\mathbf 1\not\perp a$.

Finally, in contrast to Archimedean $f$-algebras, a local pre-$d$-algebra need not be commutative.

\begin{example}\label{ex:noncomm}
Let $A$ be the space of real trigonometric polynomials on the circle $S^1=\mathbb R/2\pi\mathbb Z$ with the pointwise order, as in Example~\ref{ex:nofcomp}, fix $\theta_0\in S^1$, and define $a*b=b(\theta_0)\,a$ for $a,b\in A$.
\begin{enumerate}[(1)]
\item The multiplication is associative, since $(a*b)*c=b(\theta_0)c(\theta_0)\,a=a*(b*c)$, and products of positive elements are positive.
\item For $b\in A_+$ the operator $R_b=b(\theta_0)I$ is a Riesz homomorphism. For $a\in A_+$ and $b_1,b_2\in A$ put $M=\max(b_1(\theta_0),b_2(\theta_0))$. As in Example~\ref{ex:nofcomp}, the values $b(\theta_0)$ of upper bounds $b\in A$ of $\{b_1,b_2\}$ have infimum $M$, so $i(L_ab_1)\vee i(L_ab_2)=M\,i(a)=\inf\{i(L_ab): b\in A,\ b\ge b_1,b_2\}$. By \cite[Proposition~2.3.26(ii)]{KvG19} and \cite[Lemma~2.8]{KK22}, $L_a$ is a Riesz homomorphism. Hence $A$ is a pre-$d$-algebra.
\item The algebra $A$ is not commutative, and squares need not be positive: for $a\in A$, $a*a=a(\theta_0)a$ is not positive if $a(\theta_0)<0$ and $a(\theta)>0$ for some $\theta\in S^1$. This illustrates Theorem~\ref{thm:comm} below. As $A$ has no nonzero disjoint elements, $A$ is also a local pre-Riesz algebra. The extended multiplication on $\hat A$ is $b\cdot b'=b'(\theta_0)\,b$ $(b,b'\in\hat A)$, and $\hat A$ is a non-commutative Archimedean $d$-algebra, a function-space analogue of \cite[Example~4.1]{BH90}.
\end{enumerate}
The same construction works on every order dense subspace $A$ of $C(K)$, $K$ a compact Hausdorff space, and every point $t_0\in K$ for which evaluation at $t_0$ is a Riesz homomorphism on $A$. For instance, on $\mathbb R^3$ with the ice cream cone \cite[Section~1.7.1]{KvG19}, every extreme point $\varphi$ of the set of positive linear functionals $\psi$ on $\mathbb R^3$ with $\psi(u)=1$, where $u=(0,0,1)$, is a Riesz homomorphism \cite[Proposition~2.5.5]{KvG19}, and $a*b=\varphi(b)\,a$ $(a,b\in\mathbb R^3)$ is a pre-$d$-algebra structure on a finite-dimensional space that is not a Riesz space.
\end{example}

\section{From $d$-algebras to $f$-algebras}\label{sec:f}

Let $A$ be a pre-$d$-algebra. By Theorem~\ref{thm:ext}, its $d$-completion $\hat A$ is an Archimedean $d$-algebra. Besides the sign condition of Definition~\ref{def:pref}, we give further sufficient conditions under which $\hat A$ is an $f$-algebra, that is, under which $A$ is a pre-$f$-algebra (Theorem~\ref{thm:fcomp}). By Proposition~\ref{prop:pflocal}, each of them also implies that $A$ is local.

By Theorem~\ref{thm:BH}(2), every $d$-algebra with a unit $e>0$ is an $f$-algebra. We first show that the $d$-completion of a unital pre-$d$-algebra is an $f$-algebra.

\begin{theorem}\label{thm:unital}
Let $A$ be a pre-$d$-algebra with a multiplicative unit $e$. Then $e\in A_+$, $i(e)$ is the unit of $\hat A$, and $\hat A$ is an Archimedean $f$-algebra. Consequently $A$ is commutative.
\end{theorem}

\begin{proof}
Write $i(e)=i(e)^+-i(e)^-$ in $\hat A$. Since $\hat A$ is a $d$-algebra, Theorem~\ref{thm:BH}(1) gives, for $b\in\hat A_+$,
\[
(b\cdot i(e))^-=b\cdot i(e)^-,\qquad (b\cdot i(e))^+=b\cdot i(e)^+.
\]
For $a\in A_+$ we have $i(a)\cdot i(e)=i(ae)=i(a)\ge0$, so $i(a)\cdot i(e)^-=(i(a)\cdot i(e))^-=0$. The map $\hat A\to\hat A$, $b\mapsto b\cdot i(e)^-$, is a lattice homomorphism (as $i(e)^-\ge0$) that vanishes on $i(A_+)$, hence on $i(A)$ since $A$ is directed, hence on $\hat A$ by Lemma~\ref{lem:uniq}. Thus $b\cdot i(e)^-=0$ for all $b\in\hat A$. Consequently $b\cdot i(e)^+=b\cdot i(e)+b\cdot i(e)^-=b$ for $b\in i(A)$, and the lattice homomorphism $\hat A\to\hat A$, $b\mapsto b\cdot i(e)^+$, coincides with the identity on $i(A)$, hence on $\hat A$.

The same argument with left multiplication shows $i(e)^-\cdot b=0$ and $i(e)^+\cdot b=b$ for all $b\in\hat A$. Therefore
\[
i(e)^-=i(e)^-\cdot i(e)^+=0,
\]
where the first equality uses that right multiplication by $i(e)^+$ is the identity and the second that left multiplication by $i(e)^-$ vanishes. Hence $i(e)=i(e)^+\ge0$, and $b\cdot i(e)=i(e)\cdot b=b$ for all $b\in\hat A$, so $i(e)$ is the unit of $\hat A$. Since $i(e)\ne0$ (as $A\ne\{0\}$) we have $i(e)>0$, and $e\in A_+$ by bipositivity of $i$.

By Theorem~\ref{thm:BH}(2), a $d$-algebra with a positive unit is an $f$-algebra. Since $\hat A$ is Archimedean, it is commutative, and since $i$ is injective and multiplicative, $A$ is commutative.
\end{proof}

By Theorem~\ref{thm:BH}(2), every semiprime $d$-algebra is an $f$-algebra. The corresponding statement holds for semiprime pervasive pre-$d$-algebras.

\begin{theorem}\label{thm:semiprime}
Let $A$ be a pre-$d$-algebra that is semiprime and pervasive. Then $\hat A$ is a semiprime Archimedean $f$-algebra. Consequently $A$ is commutative.
\end{theorem}

\begin{proof}
We show that $\hat A$ is semiprime; Theorem~\ref{thm:BH}(2) then gives the $f$-algebra property, and commutativity of $A$ follows as in Theorem~\ref{thm:unital}. Let $b\in N(\hat A)$. By Theorem~\ref{thm:BH}(4), $N(\hat A)$ is an $l$-ideal, so $|b|\in N(\hat A)$. Suppose $b\ne0$. Since $A$ is pervasive, there is $a\in A$ with $0<i(a)\le|b|$. As $N(\hat A)$ is an order ideal, $i(a)\in N(\hat A)$, say $i(a)^n=0$ with $n\in\mathbb N$. Then $i(a^n)=i(a)^n=0$, so $a^n=0$ by injectivity of $i$, and $a=0$ because $A$ is semiprime. This contradicts $i(a)>0$. Hence $N(\hat A)=\{0\}$.
\end{proof}

If the pre-$d$-algebra $A$ is a Riesz space, then $\hat A=A$ and $A$ is an Archimedean $d$-algebra, since Definition~\ref{def:alg} requires $A$ to be Archimedean. Theorems~\ref{thm:unital} and~\ref{thm:semiprime} then reduce to \cite[Theorem~1.4]{BH90} for Archimedean $d$-algebras, whereas \cite[Theorem~1.4]{BH90} itself does not require the Archimedean property; in Theorem~\ref{thm:unital}, however, the positivity of the unit is a conclusion and not a hypothesis. The unit in Theorem~\ref{thm:unital} and semiprimeness in Theorem~\ref{thm:semiprime} cannot be omitted: the $d$-algebra of Example~\ref{ex:R2} is a Riesz space, hence pervasive, but it has no unit, it is not semiprime, and it is not an $f$-algebra.

A pre-$d$-algebra is commutative if and only if its squares are positive.

\begin{theorem}\label{thm:comm}
Let $A$ be a pre-$d$-algebra. The following are equivalent:
\begin{enumerate}[(1)]
\item $A$ is commutative;
\item $a^2\ge0$ for every $a\in A$;
\item $\hat A$ is an almost $f$-algebra;
\item $\hat A$ is commutative.
\end{enumerate}
\end{theorem}

\begin{proof}
(1)$\Rightarrow$(4) is Theorem~\ref{thm:ext}. (4)$\Rightarrow$(3) is \cite[Theorem~4.3]{BH90}, and (3)$\Rightarrow$(1) follows from \cite[Theorem~2.15]{BH90}, since $\hat A$ is Archimedean, and injectivity of $i$. If $\hat A$ is commutative, its squares are positive by Theorem~\ref{thm:BH}(3), applied to the Archimedean $d$-algebra $\hat A$, so $i(a^2)=i(a)^2\ge0$ and $a^2\ge0$ by bipositivity; this proves (4)$\Rightarrow$(2).

It remains to prove (2)$\Rightarrow$(3). By \cite[Proposition~1.3]{BH90}, $\hat A$ is an almost $f$-algebra if and only if $b^+b^-=b^-b^+=0$ for all $b\in\hat A$. Let first $a\in A$. Then $i(a)^2=i(a^2)\ge0$, and Theorem~\ref{thm:BH}(1) gives
\[
0=\bigl(i(a)^2\bigr)^-=i(a)^+i(a)^-+i(a)^-i(a)^+,
\]
so $i(a)^+i(a)^-=i(a)^-i(a)^+=0$, both terms being positive. Now let $b\in\hat A$ and choose $a_{jk}\in A$ with \eqref{eq:red} as in the proof of Theorem~\ref{thm:sign}. For positive factors, multiplication in $\hat A$ preserves finite suprema and infima in each variable (Theorem~\ref{thm:ext}). Hence
\[
b^+b^-=\bigvee_{j,k}\ \bigwedge_{k',j'}\ i(a_{jk'})^+\,i(a_{j'k})^-\ \le\ \bigvee_{j,k} i(a_{jk})^+\,i(a_{jk})^-=0,
\]
where the inequality takes $k'=k$ and $j'=j$ in the infimum. The same argument gives $b^-b^+=0$.
\end{proof}

Theorem~\ref{thm:comm} is the pre-Riesz version of \cite[Theorem~4.5]{BH90}; condition (2) refers only to $A$.

Finally, the $f$-algebra property of the $d$-completion can be characterised by the van Haandel extensions of the multiplication operators.

\begin{theorem}\label{thm:char}
Let $A$ be a pre-$d$-algebra. Then $\hat A$ is an $f$-algebra if and only if $\hat L_c,\hat R_c\in\Orth(\hat A)$ for every $c\in A_+$.
\end{theorem}

\begin{proof}
For $c\in A_+$, by Remark~\ref{rem:vH} the maps $\hat A\to\hat A$ given by $b\mapsto b\cdot i(c)$ and $b\mapsto i(c)\cdot b$ are $\hat R_c$ and $\hat L_c$. If $\hat A$ is an $f$-algebra, multiplication by any element is an orthomorphism, so $\hat L_c,\hat R_c\in\Orth(\hat A)$.

Conversely, assume $\hat L_c,\hat R_c\in\Orth(\hat A)$ for all $c\in A_+$, and let $b_1,b_2\in\hat A$ with $b_1\perp b_2$, $b_1\ge0$. The band $\{b_2\}^d$ is a Riesz subspace of $\hat A$, and the map $\hat A\to\hat A$, $b\mapsto b_1\cdot b$, is a lattice homomorphism. Hence
\[
M=\{b\in\hat A:\ b_1\cdot b\perp b_2\}
\]
is a Riesz subspace of $\hat A$. For $c\in A_+$ we have $b_1\cdot i(c)=\hat R_c(b_1)\perp b_2$ because $\hat R_c$ is band preserving, so $i(A_+)\subseteq M$ and, $A$ being directed, $i(A)\subseteq M$. Since $\hat A$ is the Riesz subspace generated by $i(A)$, $M=\hat A$. Thus $b_1\cdot b\perp b_2$ for all $b\in\hat A$, and symmetrically $b\cdot b_1\perp b_2$. For general $b_1\perp b_2$ apply this to $b_1^+$ and $b_1^-$. Hence left and right multiplications are orthomorphisms of $\hat A$, i.e.\ $\hat A$ is an $f$-algebra.
\end{proof}

We now apply Theorem~\ref{thm:char} to order dense subspaces of $C(K)$, where $K$ is a compact Hausdorff space. For such spaces the van Haandel extensions of the multiplication operators are weighted composition operators (Theorem~\ref{thm:vI}), and we show that local pre-Riesz algebras with sufficiently many disjoint elements are pre-$f$-algebras. Recall from \cite[Section~8]{KSvG21} that a subspace $E\subseteq C(K)$ satisfies condition (SD) if
\begin{itemize}
\item[(SD)] for all $p_1\ne p_2$ in $K$ there are $x_1,x_2\in E_+$ with $x_1\perp x_2$, $x_1(p_1)\ne0$ and $x_2(p_2)\ne0$,
\end{itemize}
that is, the points of $K$ are separated by disjoint positive elements of $E$.
Condition (SD) implies that $E$ strongly separates the points of $K$ \cite[Lemma~48]{KSvG21}. It holds for every \emph{pointwise pervasive} subspace $E$ of $C(K)$, that is, if for every $p\in K$ and every neighbourhood $U$ of $p$ there is $x\in E_+$ with $x(p)\neq0$ and $\{q\in K: x(q)\neq0\}\subseteq U$ \cite[Proposition~50]{KSvG21}; in particular, it holds for $C^k[0,1]$ (see \cite[Example~44]{KSvG21} for $k=1$). It fails for every subspace without nonzero disjoint positive elements \cite[Example~45]{KSvG21}. Under condition (SD), positive operators that are disjointness preserving on the positive cone are Riesz* homomorphisms.

\begin{theorem}[{\cite[Theorem~49]{KSvG21}}]\label{thm:KSvG}
Let $K$ be a compact Hausdorff space and let $E\subseteq C(K)$ be order dense and satisfy \textup{(SD)}. Then every positive linear map $T\colon E\to E$ that is disjointness preserving on $E_+$ is a Riesz* homomorphism.
\end{theorem}

For order dense subspaces of $C(K)$ that strongly separate the points of $K$, the $f$-algebra property of the $d$-completion is equivalent to a pointwise condition.

\begin{proposition}\label{prop:pw}
Let $K$ be a compact Hausdorff space, let $A\subseteq C(K)$ be an order dense subspace that strongly separates the points of $K$, and let $A$ be a pre-$d$-algebra. Then $\hat A$ is an $f$-algebra if and only if, for all $a\in A$, $c\in A_+$ and $q\in K$,
\[
a(q)=0\ \Longrightarrow\ (ca)(q)=(ac)(q)=0.
\]
\end{proposition}

\begin{proof}
Let $c\in A_+$. Since $L_c$ is a Riesz* homomorphism (Remark~\ref{rem:vH}), Theorem~\ref{thm:vI} and the subsequent paragraph give $w_c\in C(K)$, $w_c\ge0$, and $\alpha_c\colon K\to K$ with $(ca)(q)=w_c(q)\,a(\alpha_c(q))$ for $a\in A$, $q\in K$, and $\hat L_cb=w_c\cdot(b\circ\alpha_c)$ for $b\in\hat A$.

Assume the pointwise condition, and let $q\in K$ with $w_c(q)>0$. If $\alpha_c(q)\ne q$, choose $b\in A$ with $b(q)=0$ and $b(\alpha_c(q))=1$; then $(cb)(q)=w_c(q)\ne0$, a contradiction. Hence $\alpha_c(q)=q$ wherever $w_c(q)>0$, so $\hat L_cb=w_cb$ for $b\in\hat A$. As the lattice operations of $\hat A$ are pointwise, $\hat L_c$ is band preserving, i.e.\ $\hat L_c\in\Orth(\hat A)$. The same holds for $\hat R_c$, and Theorem~\ref{thm:char} shows that $\hat A$ is an $f$-algebra.

Conversely, let $\hat A$ be an $f$-algebra, and suppose $a(q)=0$ but $(ca)(q)\neq0$ for some $a\in A$, $c\in A_+$, $q\in K$. Then $w_c(q)>0$ and $p:=\alpha_c(q)\neq q$. Choose $b\in A$ with $b(q)=0$ and $b(p)=1$, and $b'\in A$ with $b'\ge\mathbf 1$, which exists since $A$ is order dense and hence majorizing in $C(K)$. For $0<\varepsilon<1/b'(p)$, $b_1=(b-\varepsilon b')^+$ and $b_2=(b-\varepsilon b')^-$ are disjoint elements of $\hat A$ with $b_1(p)>0$ and $b_2(q)>0$. Then $(\hat L_cb_1)(q)=w_c(q)b_1(p)>0$, so $\hat L_cb_1$ is not disjoint from $b_2$, contradicting $\hat L_c\in\Orth(\hat A)$. The case of $ac$ is the same with $\hat R_c$.
\end{proof}

The pointwise condition is the multiplicative form of the strong locality condition (5.1) in \cite[Remark~5.1.5]{KvG19}. In Example~\ref{ex:nofcomp}, $(c*a)(q)=c(\theta_0)a(\theta_0)$ for all $a,c\in A$ and $q\in S^1$, so the condition fails for $c=\mathbf 1$ whenever $a(\theta_0)\ne0=a(q)$.

We now show that under condition (SD) locality and the pre-$f$ property coincide.

\begin{theorem}\label{thm:main}
Let $K$ be a compact Hausdorff space and let $A$ be a pre-Riesz algebra that is an order dense subspace of $C(K)$ satisfying \textup{(SD)}. The following are equivalent.
\begin{enumerate}[(i)]
\item $A$ is local.
\item $A$ is a pre-$f$-algebra.
\end{enumerate}
In this case $A$ is commutative, and for every $c\in A_+$ there is $w_c\in C(K)_+$ with
\[
(ca)(q)=(ac)(q)=w_c(q)\,a(q)\qquad(a\in A,\ q\in K).
\]
\end{theorem}

\begin{proof}
(ii)$\Rightarrow$(i) is Proposition~\ref{prop:pflocal}.

(i)$\Rightarrow$(ii). Let $c\in A_+$. The operator $L_c$ is positive and local, hence disjointness preserving, in particular on $A_+$. By Theorem~\ref{thm:KSvG} it is a Riesz* homomorphism, and since \textup{(SD)} implies strong separation of points, Theorem~\ref{thm:vI} gives $w_c\in C(K)$, $w_c\ge0$, and $\alpha_c\colon K\to K$ with $(ca)(q)=w_c(q)\,a(\alpha_c(q))$ for $a\in A$, $q\in K$.

We show that $\alpha_c(q)=q$ whenever $w_c(q)>0$. Suppose $w_c(q)>0$ and $p:=\alpha_c(q)\neq q$. By \textup{(SD)} there are $a_1,a_2\in A_+$ with $a_1\perp a_2$, $a_1(p)\ne0$ and $a_2(q)\ne0$. Locality gives $ca_1\perp a_2$, i.e.\ $(ca_1)(t)\,a_2(t)=0$ for all $t\in K$. But $(ca_1)(q)\,a_2(q)=w_c(q)\,a_1(p)\,a_2(q)\neq0$, a contradiction. Hence $(ca)(q)=w_c(q)a(q)$ for all $a\in A$ and $q\in K$, that is, $L_c$ is the restriction to $A$ of the multiplication operator $M_{w_c}\colon x\mapsto w_cx$ on $C(K)$.

The multiplication operator $M_{w_c}$ is a positive orthomorphism of $C(K)$, hence an order continuous lattice homomorphism, and it maps $A$ into $A$. By Lemma~\ref{lem:RH}, applied to the cover $C(K)$, $L_c$ is a Riesz homomorphism of $A$. The same argument applies to $R_c$ and gives $v_c\in C(K)$, $v_c\ge0$, with $(ac)(q)=v_c(q)a(q)$ for $a\in A$, $q\in K$. Hence $A$ is a pre-$d$-algebra, and it satisfies the pointwise condition of Proposition~\ref{prop:pw}, since $(ca)(q)=w_c(q)a(q)$ and $(ac)(q)=v_c(q)a(q)$ vanish whenever $a(q)=0$. By Proposition~\ref{prop:pw}, $\hat A$ is an $f$-algebra, so $A$ is a pre-$f$-algebra by Theorem~\ref{thm:fcomp}. By the same theorem $A$ is commutative, so $(ac)(q)=(ca)(q)=w_c(q)\,a(q)$.
\end{proof}

The following example concerns subalgebras of $C[0,1]$ with the pointwise product and illustrates Theorems~\ref{thm:unital} and~\ref{thm:main}.

\begin{example}\label{ex:fun}
Let $A\subseteq C[0,1]$ be a linear subspace that contains the polynomials and is closed under pointwise multiplication, for instance the polynomials, $C^k[0,1]$ or the real-analytic functions $C^\omega[0,1]$, with the pointwise order and the pointwise product. By the Weierstrass theorem, $A$ is order dense in $C[0,1]$; hence $A$ is an Archimedean pre-Riesz space, and $\hat A$ is the Riesz subspace of $C[0,1]$ generated by $A$. For $c\in A_+$, multiplication by $c$ is an order continuous lattice homomorphism of $C[0,1]$ that maps $A$ into $A$; by Lemma~\ref{lem:RH}, applied to the cover $C[0,1]$, $L_c=R_c$ is a Riesz homomorphism of $A$. So $A$ is a pre-$d$-algebra with unit $\mathbf 1$, and Theorem~\ref{thm:unital} shows that $\hat A$ is an $f$-algebra. The multiplication of $\hat A$ is the pointwise product: applying Lemma~\ref{lem:uniq} with $F=C[0,1]$, first to the maps $\hat A\to C[0,1]$ given by $b\mapsto b\cdot c$ and $b\mapsto bc$ for $c\in A_+$, and then to the maps $\hat A\to C[0,1]$ given by $b\mapsto b_1\cdot b$ and $b\mapsto b_1b$ for $b_1\in\hat A_+$, we obtain $b_1\cdot b_2=b_1b_2$ for all $b_1,b_2\in\hat A$. In particular, the Riesz subspace of $C[0,1]$ generated by $A$ is closed under pointwise multiplication. The spaces listed above are not Riesz spaces. For $A=C^k[0,1]$ the $f$-algebra property of $\hat A$ also follows from Theorem~\ref{thm:main}, with $w_c=c$: $A$ satisfies (SD), and $A$ is local, since for $a,b\in A$ and $c\in A_+$, $a\perp b$ implies $a(q)b(q)=0$ and hence $(ca)(q)b(q)=0$ for all $q\in[0,1]$. The polynomials and $C^\omega[0,1]$ have no nonzero disjoint elements, so they do not satisfy (SD).
\end{example}

\end{document}